\documentclass[12pt]{amsart}
\usepackage{amssymb,amsmath}

\def\Om{\Omega}

\numberwithin{equation}{section}
\def\db{\bar\partial}
\def\db*{\bar\partial^*}
\def\T{\text}
\def\simgeq{\underset\sim>}
\def\simleq{\underset\sim<}
\def\1#1{\overline{#1}}
\def\2#1{\widetilde{#1}}
\def\3#1{\widehat{#1}}
\def\4#1{\mathbb{#1}}

\def\5#1{\frak{#1}}
\def\6#1{{\mathcal{#1}}}
\def\C{{\4C}}
\def\R{{\4R}}

\def\Z{{\4Z}}

\def\di{\partial}
\def\dib{\bar\partial}
\begin{document}
\title[Propagation at the vertex of a sector]{Propagation at the vertex of a sector}
\author[S.~Pinton ]{Stefano Pinton }
\email{pinton@math.unipd.it}
\maketitle
\begin{abstract}
We discuss holomorphic extension across a boundary point in terms of sector property. The point is of infinite type and the sector is accordingly ``cusped" at the vertex.

\noindent
32F10, 32F20, 32N15, 32T25 
\end{abstract}
\def\Giialpha{\mathcal G^{i,i\alpha}}
\def\cn{{\C^n}}
\def\cnn{{\C^{n'}}}
\def\ocn{\2{\C^n}}
\def\ocnn{\2{\C^{n'}}}
\def\const{{\rm const}}
\def\rk{{\rm rank\,}}
\def\id{{\sf id}}
\def\aut{{\sf aut}}
\def\Aut{{\sf Aut}}
\def\CR{{\rm CR}}
\def\GL{{\sf GL}}
\def\Re{{\sf Re}\,}
\def\Im{{\sf Im}\,}
\def\codim{{\rm codim}}
\def\crd{\dim_{{\rm CR}}}
\def\crc{{\rm codim_{CR}}}
\def\phi{\varphi}
\def\eps{\varepsilon}
\def\d{\partial}
\def\a{\alpha}
\def\b{\beta}
\def\g{\gamma}
\def\G{\Gamma}
\def\D{\Delta}
\def\Om{\Omega}
\def\k{\kappa}
\def\l{\lambda}
\def\L{\mathcal L}
\def\z{{\bar z}}
\def\w{{\bar w}}
\def\Z{{\1Z}}
\def\t{{\tau}}
\def\th{\theta}
\emergencystretch15pt
\frenchspacing
\newtheorem{Thm}{Theorem}[section]
\newtheorem{Cor}[Thm]{Corollary}
\newtheorem{Pro}[Thm]{Proposition}
\newtheorem{Lem}[Thm]{Lemma}
\theoremstyle{definition}\newtheorem{Def}[Thm]{Definition}
\theoremstyle{remark}
\newtheorem{Rem}[Thm]{Remark}
\newtheorem{Exa}[Thm]{Example}
\newtheorem{Exs}[Thm]{Examples}
\def\Label#1{\label{#1}}
\def\bl{\begin{Lem}}
\def\el{\end{Lem}}
\def\bp{\begin{Pro}}
\def\ep{\end{Pro}}
\def\bt{\begin{Thm}}
\def\et{\end{Thm}}
\def\bc{\begin{Cor}}
\def\ec{\end{Cor}}
\def\bd{\begin{Def}}
\def\ed{\end{Def}}
\def\br{\begin{Rem}}
\def\er{\end{Rem}}
\def\be{\begin{Exa}}
\def\ee{\end{Exa}}
\def\bpf{\begin{proof}}
\def\epf{\end{proof}}
\def\ben{\begin{enumerate}}
\def\een{\end{enumerate}}
\def\dotgamma{\Gamma}
\def\dothatgamma{ {\hat\Gamma}}

\def\simto{\overset\sim\to\to}
\def\1alpha{[\frac1\alpha]}
\def\T{\text}
\def\R{{\Bbb R}}
\def\I{{\Bbb I}}
\def\C{{\Bbb C}}
\def\Z{{\Bbb Z}}
\def\Fialpha{{\mathcal F^{i,\alpha}}}
\def\Fiialpha{{\mathcal F^{i,i\alpha}}}
\def\Figamma{{\mathcal F^{i,\gamma}}}
\def\Real{\Re}
%
%
%
\section{Introduction}
\Label{s0}
This paper looks at the holomorphic extension across a smooth real hypersurface, the boundary $b\Om$ of a domain of the complex space, about a point that we fix as the origin $0$. Our program is to give the analogous of the result of Baouendi and Treves \cite{BT84} and Baracco, Zaitsev and Zampieri \cite{BZZ07} on holomorphic estension across the boundary in terms of the ``sector property". In our setting the boundary may have infinite type and the sector, described by the inverse function of the type, is accordingly singular at its vertex $0$. We also prefere to use the language of propagation  instead of forced extension; however, the two points of view do not differ substantially. In coordinates $(z,z',w)\in \C\times \C^{n-2}\times \C$, $z=x+iy,\,\,z'=x'+iy',\,\,w=r+is$, let  $b\Om$ be defined by $s=h$ for $h(0)=0$, $dh(0)=0$ and, for a holomorphic function $F(z),\,\,z\in\C$, suppose that $h$ satisfies $h|_{z'=0}=O(F(|z|))$ (with some minor additional requirements such as \eqref{1.1,6}, \eqref{1.1,7} and \eqref{increment} below). The sector property of $b\Om$ over $S_\alpha$, parametrized by $ F^*_\alpha:=F^*(\epsilon z^\alpha),\,\,z=1-\tau\in1-\Delta$ (where $F^*$ is the inverse to $F$ and $\Delta$ is the standard disc) consists in $h|_{S_\alpha}\le0$ for $\alpha>1$. Our result is that, 
\begin{itemize}
\item[(i)] 
If the sector property holds, then the disc $\mathcal S_\alpha$ attached to $b\Om$ over $S_\alpha$ is a propagator of holomorphic extendibility from $\Om$ to $s<h$ at the vertex $0$.
\item[(ii)] Under some additional condition (cf. \eqref{2.2} below), the sector property is  necessary for holomorphic extendibility. If it is not satisfied, there is a fundamental system of neighborhoods of  $\Om\cup\mathcal S_\alpha$ which are pseudoconvex. Thus propagation cannot occur.
\end{itemize}
Our result is best understood in the model situation in which $n=2$. We first treat the finite type.
\be
\Label{e0.1}
({\it Finite  type.}) We consider the domain $\Om\subset\C^2$ defined by $s>h$ for $h=|z|^{2m}+c|z|^{2m-2p}(\Re z)^{2p}$
and discuss the choice of $c$ which yields the sector property.
For this, we set $F^*:=z^{\frac1{2m}}$ and $S_\alpha=F^*_\alpha(1-\Delta)$. Then $h<0$ in a sector $S_\alpha$ for $\alpha>1$ if and only if $c>\cos^{-1}(\frac{2p\pi}{4m})$ (cf \cite{BZ08} Proposition 4.1). Under this assumption,  $S_\alpha\times\{0\}$ is contained in $\bar\Om$ and thus it is a propagator  by Theorem~\ref{t1.1} and Remark~\ref{r1.1} below (or also  by \cite{BT84}). When $c\le \cos^{-1}(\frac{2p\pi}{4m})$ (and in case $p$ is a divisor of $m$) there are holomorphic functions on $\Om$ which do not extend down at $0$ (cf. \cite{BZ08} Corollary 4.3) since, in new holomorphic coordinates, it is contained in the half space $s>0$. Notice that the condition on $c$ for the pseudoconvexity of $\Om$ is different; for instance for $m=2$ and $p=1$, this is $c\le\frac43$ whereas the condition for the sector property is $c\le\cos^{-1}(\frac\pi4)=\sqrt2$. This means that for the intermediate values $\frac43<c<\sqrt2$ we have holomorphic extension at points arbitrarily close to $0$ but not at $0$.
\ee
We pass to discuss the infinite type. We consider pairs of functions such as $(F,F^*)=(e^{-\frac 1{z^a}},\,\frac1{(-\log\,z)^{\frac1a}})$ or $(F,F^*)=(e^{-e^{\frac1{z^a}}},\,\frac1{\log(-\log\,z)^{\frac1a}})$.
\bl
\Label{l2.2}
 The sector $S_\alpha$ parametrized by $F^*_\alpha=F^*(\epsilon(1-\tau)^\alpha)$ in the two respective cases listed above, satisfies
\begin{equation}
\Label{nova}
\begin{cases}
S_\alpha\sim\{z\in \C^+:\,\,|y|<\alpha x^{a+1}\},
\\
\T{ (resp. $S_\alpha\sim\{z\in\C^+:\,\,|y|<c_{\alpha,\,x}x^{a+1}e^{-\frac1{x^a}}e^{-e^{\frac1{|x|^a}}}\}$ for $c_{\alpha,\,x}=1+\frac{\log\alpha}{\log(-\log x)}$.)}
\end{cases}
\end{equation}
\el
\bpf
In proving \eqref{nova}  we neglect the factor $\epsilon$ which is irrelevant at the vertex. We start from the first of \eqref{nova} and observe that $S_\alpha$ at $z=0$ is the $\frac1{\alpha^{\frac1a}}$-homotetic set of $S_1=\{\frac1{(-\log(1-\tau))^{\frac1a}}:\,\,\tau\in\Delta\}$ at $\tau=1$. Now, to describe the points $z=x+iy\in b S_1$, we write $\tau=e^{i\theta}$, for $|\theta|$  small, and get
\begin{equation*}
\begin{split}
x+iy&\sim \frac{1}{(-\log (|1-e^{i\theta}|))^{\frac1a}}+i\frac{\arg (1-e^{i\theta})}{(-\log (|1-e^{i\theta}|))^{\frac1a+1}}
\\&\sim\frac1{(-\log |\theta| )^{\frac1a}}+i\frac{\T{arctg}\,(\frac2\theta)}{(-\log |\theta| )^{\frac1a+1}}
\\
&\sim \frac1{(-\log |\theta| )^{\frac1a}}\pm i\frac1{(-\log |\theta| )^{\frac1a+1}}.
\end{split}
\end{equation*}
Thus $z=x+iy\in bS_\alpha$ if and only if $|y|\sim x^{\a+1}$. Taking the $\frac1{\alpha^a}$-homotetic set, the inequality changes into $|y|\sim\alpha x^{\a+1}$.

To prove the second of \eqref{nova}, an easy calculation shows that $S_\alpha$ is approximately the $\frac1{c_{\alpha,\,x}^{\frac1a}}$-homotetic of $S_1$ for $c_{\alpha,\,x}=1+\frac{\log\alpha}{\log(-\log x)}$.
We observe that, since $-\log(1-e^{i\theta})\sim-log\theta^2+i\theta$, then $\arg(-\log(1-e^{i\theta}))\sim\frac\theta{-\log\theta^2}$. Thus,   the points $z=x+iy\in b S_1$ are described by
\begin{equation*}
\begin{split}
x+iy&=\frac1{(\log(-\log|\theta| ))^{\frac1a}}+i\frac{\arg(-\log(1-e^{i\theta})}{(\log(-\log|\theta| ))^{\frac1a+1}}
\\
&\sim \frac1{(\log(-\log|\theta| ))^{\frac1a}}+i\frac\theta{-\log|\theta| (\log(-\log|\theta| ))^{\frac1a+1}}.
\end{split}
\end{equation*}
Thus $ S_1$ is described by $|y|\simleq x^{a+1}e^{-\frac1{|x|^a}}e^{-e^{\frac1{|x|^a}}}$ and its $\frac1{c_{\alpha,\,x}^{\frac1a}}$-homotetic $S_\alpha$ by $ |y|\simleq c_{\alpha,\,x}x^{a+1}e^{-\frac1{|x|^a}}e^{-e^{\frac1{|x|^a}}}$. This completes the proof of \eqref{nova}.

\epf
\be
\Label{e0.2} {\it (Infinite type.)} If $S_\alpha$ is the sector parametrized by $F^*_\alpha$ for $\alpha>1$ (in the two respective cases),  we introduce the cut-off $\chi=\chi\Big(\frac y{2\alpha x^{a+1}}\Big)$ (resp. $\chi=\chi\Big(\frac y{2c_{\alpha,x}x^{a+1}e^{-\frac1{x^a}e^{-e^{\frac1{|x|^a}}}}}\Big)$)
which has support in $S_{2\alpha}$ and is $1$ on $S_\alpha$.  We choose $b$ such that $a<b<a(a+1)$ and
 consider the domain $\Om$ of infinite type defined by $s>h$ for the two choices
\begin{equation*}
\begin{cases}
h=e^{-\frac1{|y|^a}}-\chi\Big(\frac y{2\alpha x^{a+1}}\Big)\Re e^{-\frac 1{z^b}},
\\
(\T{ resp. $h=e^{-e^{\frac1{|y|^a}}}-\chi\Big(\frac y{2c_{\alpha,x}x^{a+1}e^{-\frac1{x^a}e^{-e^{\frac1{|x|^a}}}}}\Big)\Re e^{-e^{\frac 1{z^b}}})$}.
\end{cases}
\end{equation*}
We observe that 
$\Re e^{-\frac 1{z^b}}$ is positive on supp $\chi$ and the same is true for the double exponential. Moreover, 
because of $b<(a+1)a$, we have for any $\alpha>1$ and in a suitable neighborhood of $0$
\begin{equation*}
\begin{cases}
\Re e^{-\frac 1{z^b}}>e^{-\frac1{|z|^a}}\quad\T{ on $S_\alpha$},
\\
\T{(resp. }\Re e^{-e^{\frac 1{z^b}}}>e^{-e^{\frac1{|z|^a}}}\quad\T{on } S_\alpha).
\end{cases}
\end{equation*}
Thus $\Om\supset S_\alpha\times\{0\}$. Also, because of $b>a$, we have satisfied \eqref{1.1} and \eqref{1.1,6} (whereas \eqref{1.1,7} and \eqref{increment} are obvious). Hence the hypotheses of Theorem~\ref{t1.1} are all fulfilled  and $S_\alpha\times\{0\}$ is a propagator of holomorphic extendibility at $0$; in fact, a stronger property holds, that is, forced extension.

 Note that, when $a\ge1$, it has been proved in \cite{BKZ12} that  the half line $\R^+$ itself is a propagator. We can also observe that $e^{-e^{\frac1{y^a}}}<e^{-\frac 1{y^c}}$ for any $c$ and thus in particular for $c>1$; thus $\R^+$ is also a propagator for the domain  $s>e^{-e^{\frac1{y^a}}}$ regardless who is $a$. However, this requires different tools from the present context.
\ee


\be
{\it (Failure of the sector property in the exponentially degenerate type.)}
Let  us consider the  ``tube" domain defined by $s>e^{-\frac1{|y|^a}}$ for  $a<1$; is $\R^+$ still  a propagator? To decide it,  we take $b$ such that  $(b+1)a<b$, consider the holomorphic function $e^{-\frac 1{z^b}}$, its inverse $\frac1{(-\log)^{\frac1b}}$ and the related sector
$
\{z\in\C^+:\,\,|y|< x^{b+1}\}.
$
(Here the use of $\alpha$ has become irrelevant.)
Note that the choice of the holomorphic function is no longer related to the vanishing order of the boundary.
We introduce the modified domain $ \Om$ defined by $s>h$ for $ h=e^{-\frac1{|y|^a}}-\chi\Big(\frac y{2x^{b+1}}\Big)\Re e^{-\frac 1{z^b}}$,
where $\chi\Big(\frac y{2x^{b+1}}\Big)$ is $1$ for $y<x^{b+1}$ and has support in $y<2x^{b+1}$. 
We have $\di_z\di_{\bar z} h>\frac{|\di_z e^{-\frac1{z^b}}|^2}{|e^{-\frac1{z^b}}|}$ on supp $\dot\chi$ and therefore, by Proposition~\ref{p2.1} below,  $h$ is plurisubharmonic, $\Om$ is pseudoconvex, and $ S_\alpha\times\{0\}$ is not a propagator for $ \Om$ since we can subtract to $h$ a further term $\chi\Re e^{-\frac 1{z^b}}$ without destroying its pseudoconvexity.  In particular, $\R^+$ is not a propagator of extendibility for holomorphic functions on the initial tube domain. (To prove directly the plurisubharmonicity of $ h$, it suffices to notice that, over supp $\dot\chi$, 
\begin{equation*}
\begin{split}
 \di_z\di_{\bar z} e^{-\frac1{|y|^a}}&\simgeq \frac{e^{-\frac1{x^{a(b+1)}}}}{x^{(b+1)(2a+2)}}
\\
&\underset{\T{since $a(b+1)<b$}}\simgeq 
\frac1{x^{b+1}}e^{-\frac c{x^b}}\quad\T{for $c<1$}
\\
&\simgeq \di_z\di_{\bar z}\Big(\chi\Big(\frac y{2x^{b+1}}\Big)\Re e^{-\frac1{z^b}}\Big).
\end{split}
\end{equation*}

The conclusion does  not contradict Theorem~\ref{t1.1}, because, in the present case, the sector is not calibrated on the type of $b\Om$ by our choice of $b>a$. If, instead, 
$a<b<(a+1)a$, the sector $\{z\in\C^+:\,y<x^{b+1}\}$ is a propagator for the domain defined as in Example~\ref{e0.2}; more precisely,  we have forced extension.

\ee                                                                                                                                                                                                                                                                                                                                                                                                                                                                                                                                                                                                                                                                                                                                               

\section{Propagation of holomorphic extendibility}
\Label{s1}
Let $F^*$ be a holomorphic, injective function defined on the points $\epsilon z^\alpha$, for $z\in \C^+:=\{z\in C:\,\Re z>0\}$ whose range contains a neighborhood of $0$ in the half line $\R^+$.
For positive $\alpha$ and for $F_\alpha^*:=F^*(\epsilon z^\alpha)$, let $S=S_\alpha$ be the sector $\{z=F^*_\alpha(1-\tau),\,\,\tau\in\Delta\}$ parametrized by $F^*_\alpha$ over the standard disc $\Delta$; this  is possibly singular  at $\tau=1$.
 We denote by $F$ the inverse to $F^*$. The pairs of functions $(F^*,F)$ that we have in mind are $(z^{\frac1{2m}},z^{2m})$ or else $(\frac1{(-\log z)^{\frac1a}},\,e^{-\frac1{z^a}})$.   Let $\Om$ be
 a domain in $\C^n$ defined by $s>h$ for $h\le cF(|z|)$  when $z$ moves in a complex curve. We take coordinates $(z,z',w)\in \C\times \C^{n-2}\times \C,\,\,z=x+iy, \,\,z'=x'+iy',\,\,w=r+is$, assume that $b\Om$ is defined by $s=h(z,z',r)$ for $h(0)=0$, $dh(0)=0$ and that the complex curve is the $z$-axis. With the notation $h_{r }:=h(\cdot,0,r )$,  we are assuming  that
\begin{equation}
\Label{1.1}
h_r|_{z'=0}=O(F(|z|)\quad\T{uniformly in $r$.}
\end{equation}
In particular,
\begin{equation}
\Label{1.1,5}
|h_r||_{S_\alpha\times\{0\}}\le c|\theta|^\alpha,\,\,\T{ for $e^{i\theta}\in\di\Delta$ uniformly in $r$}.
\end{equation}
Let $D^1$ and $D^2$ denote various first and second order derivatives; we assume, together with \eqref{1.1} and uniformly in $r$
\begin{equation}
\Label{1.1,6}
D^1h_r|_{S_\alpha\times\{0\}}=O(|\di_zF|),\,\,D^2h_r|_{S_\alpha\times\{0\}}=O(|\di_z^2F|).
\end{equation} 
We also assume
\begin{equation}
\Label{1.1,7}
|\di^2_zF^*|\simleq \frac{|\di_zF^*|}{|z|},\quad |\di^2_zF|\simleq \frac{|\di_zF|}{|z|}.
\end{equation}
Again, \eqref{1.1}, \eqref{1.1,6} and \eqref{1.1,7} hold for our two main models.
We have
\bl
\Label{l1.1}
Let $1<\alpha<2$ and set $\beta:=\alpha-1$; then
$h_{r }(F^*_\alpha(1-e^{i\theta}))\in C^{1,\beta}$.
\el
\bpf
Since $F^*_\alpha(1-e^{i\theta})$ is singular only at $\theta=0$, by the Hardy-Littlewood lemma it suffices to prove that $\frac d{d\theta}h_r(F^*_\alpha)=O(\theta^{\alpha-1})$ and  $\frac{d^2}{d\theta^2 }h_r(F^*_\alpha)=O(\theta^{\alpha-2})$ at $\theta=0$. Now, for $z=1-e^{i\theta}$, we have
\begin{equation*}
\begin{cases}
\begin{split}
\Big|\frac{d}{d\theta}\Big(h_{r }(F^*_\alpha)\Big)\Big|&\simleq \Big|D^1h_{r }(F^*_\alpha)D^1F^*_\alpha\Big|
\\
&\simleq \Big|\di_zF(F^*_\alpha)\Big|\Big|(\di_zF^*)_\alpha\theta^{\alpha-1}\Big|
\\
&=O(|\theta|^{-1+\alpha}),
\end{split}
\\
\begin{split}
\Big|\frac{d^2}{d\theta^2}\Big(h_{r }(F^*_\alpha)\Big)|&\simleq \Big|D^1h_{r }(F^*_\alpha)D^2F^*_\alpha\Big|+\Big|D^2h_{r }(F^*_\alpha)(D^1F^*_\alpha)^2\Big|
\\
&\simleq \Big|\di_zF(F^*_\alpha)\Big|\Big|(\di^2_zF^*)_\alpha\theta^{2\alpha-2}\Big|+\Big|(\di_zF^*)_\alpha\theta^{\alpha-2}\Big|+\Big|\di_z^2F(F^*_\alpha)\Big|\Big|(\di_z F^*)_\alpha^2\theta^{2\alpha-2}\Big|
\\
&=O(|\theta|^{-2+\alpha}).
\end{split}
\end{cases}
\end{equation*}

\epf
We consider the (Bishop) equation
\begin{equation}
\Label{bishop}
u-T_1h(F^*_\alpha,0,u)=0,
\end{equation}
where $T_1$ is the Hilbert transform normalized by the value $0$  at $\tau=1$. Recall that $T_1$ is continuous on $C^{1,\beta}(b\Delta)$; thus the mapping $(\epsilon,u)\mapsto u-T_1h(F^*_\alpha,0,u),\,\,\R\times C^{1,\beta}\to C^{1,\beta}$ (for $F^*_\alpha=F^*(\epsilon z^\alpha)$) is differentiable  and its differential in $u$ is $1-T_1\di_zh\sim1$. By the implicit function theorem, we have that, for $\epsilon$ small, there is a unique solution $u=u(e^{i\theta})$ to \eqref{bishop}. Moreover, if we consider a 1-parameter family of deformations $h_\eta$ of $h$  so that $\eta\mapsto h_\eta(F^*_{\alpha},0,r),\,\,\R\to C^{1,\beta}(b\Delta)$ is $C^k$ uniformly with respect to $r $, then $\eta\mapsto u_\eta$ is also $C^k$. 
We put $v|_{b\Delta}:=T_1u|_{b\Delta}$ and use the same notation $u,\,v$ for the harmonic extensions from $b\Delta$ to $\Delta$. 
We define $\mathcal S=\mathcal S_\alpha$ by
\begin{equation}
\Label{1.3}
\mathcal S:=\{(F^*_\alpha(1-\tau),0,u(\tau)+iv(\tau)):\,\,\tau\in\Delta\}.
\end{equation}
This is the holomorphic disc attached to $ b\Om$ over the sector $S\times\{0\}$.
We assume $F(x)$, $F^*(x)$ increasing, and $\di_xF^{*}(x)$ decreasing.
We make an additional assumption. For this, we write $z=\rho e^{i\psi}$ or else $z=\sigma(1-e^{i\theta})$ for $(\rho,\psi)\in (0,\epsilon)\times(-\frac\pi2,\frac\pi2)$ or $(\sigma,\theta)\in (1-\epsilon,1)\times((\epsilon,-\epsilon)\setminus\{0\})$. They are related by the change
\begin{equation*}
\begin{cases}
\rho=2\sigma\sin\frac\theta2
\\
\psi=\frac{\pi-\theta}2.
\end{cases}
\end{equation*}
In particular
$$
\di_\sigma=2\sin(\frac\pi2-\psi)\di_\rho.
$$
With these preliminaries we suppose that  $\di_\sigma\Re F^*_\alpha(\sigma(1-e^{i\theta}))$ and $\di_\sigma\Im F^*_\alpha(\sigma(1-e^{i\theta}))$ have the properties that
\begin{equation}
\Label{increment}
\begin{cases}
\T{ (i) they keep the same sign for fixed $\theta$}
\\
\T{(ii) their absolute value is decreasing with respect to $\sigma$}.
\end{cases}
\end{equation}
It is readily seen that $F^*=z^{\frac1{2m}}$ and $F^*=\frac1{(-\log\,z)^{\frac1a}}$ satisfy \eqref{increment}.
\bt
\Label{t1.1}
In the above situation, in particular under \eqref{1.1} and \eqref{increment} and for $\alpha>1$, we further assume  that $\mathcal S$ is tangent to $ b\Om$ at the vertex, that is,
\begin{equation}
\Label{1.4}
\di_t v=0\quad\T{at $\tau=1$}.
\end{equation}
Then there is propagation of holomorphic extendibility across $b\Om$ from any point of $b\mathcal S$ to the vertex $0$.
\et
\bpf
Let holomorphic extendibility occur in the $\eta_o$-neighborhood $V_{\eta_o}$  of $\tau=-1$; we show that it also occurs at $\tau=1$. Let $\chi=\chi(1-\tau)$ be a cut-off in $V_{\eta_o}$ and let $u=u_\eta$ be the harmonic extension of the solution of the equation
$$
u-T_1(h(F^*_\alpha,0,u)-\eta\chi)=0.
$$
Then, $\eta\mapsto u_\eta+iv_\eta,\,\,\R\to C^{1,\beta}$ is $C^\infty$ since $\eta\mapsto h(F^*_\alpha,\cdot)-\chi\eta$ is also $C^\infty$. 
In the coordinate $\tau=te^{i\theta}\in\Delta$ we have therefore the Taylor expansion
 $$
 \di_tv_\eta=\di_tv\Big|_{\eta=0}+\eta\di_t\di_\eta v_\eta\Big|_{\eta=0}+o(\eta).
 $$
 Now, the first term in the right side is $0$ by \eqref{1.4} at $\tau=1$. On the other hand, the radial derivative of $\di_\eta v_\eta$ at $\tau=1$ can be calculated by means of the convergent integral
 \begin{equation*}
 \begin{split}
 \di_t\di_\eta v_\eta&\sim\di_t\int \di_\eta v_\eta \cdot\frac{1-t^2}{1+t-2\cos\theta}d\theta
\\ 
 &\sim-\di_t\int \chi\cdot\frac{1-t^2}{1+t-2\cos\theta}d\theta
 \\
 &\sim+\int\frac\chi{1-\cos\theta}d\theta=c,\quad \T{at $\tau=1$}.
 \end{split}
 \end{equation*}
 Thus, for $\eta$ small, $\di_tv_\eta>0$ at $\tau=1$; we pick up such $v=v_\eta$. After reparametrization $z=F^*_\alpha(1-\tau)$, and with $F$ denoting as always the inverse to $F^*$,  our temporary conclusion is that $v_\eta$ satisfies
 \begin{equation}
 \Label{1.5}
 v_\eta(t)\le -cF(t)\quad t\in(1-\epsilon,1).
 \end{equation}
 We are now tempted to move the vertex of the sector $S$ to $-\epsilon$ and to attach to $b\Om$ an $\epsilon$-parameter family of discs over the sectors $-\epsilon +S$. But, over the new vertex $-\epsilon$, we do not have any more the condition \eqref{1.1} 
 in which the singularity of $F^*$ is balanced by the vanishing order of the defining function $h$ of $b\Om$. Instead, we use approximation by the smooth sectors $S_\nu$ parametrized by
 $$
 F^*_{\alpha\,\nu}=F^*_\alpha(1-\tau+\frac1\nu)-F^*_\alpha(\frac1\nu),\quad \tau\in \Delta.
 $$
We have
$$
h\Big(F^*_\alpha(1-e^{i\theta}),r \Big)\underset{\T{\eqref{1.1}}}\simleq |\theta|^{\alpha}.
$$
We also have
\begin{equation}
\begin{split}
\Big|F^*_\alpha(1-e^{i\theta}+\frac1\nu)-F^*_\alpha(\frac1\nu)\Big| &= \Big|\int_0^1\di_\sigma F^*_\alpha(\sigma(1-e^{i\theta})+\frac1\nu)d\sigma\Big|
\\
&\underset{\T{\eqref{increment} (i)}}=\int_0^1\Big|\di_\sigma F^*_\alpha(\sigma(1-e^{i\theta})+\frac1\nu)\Big|d\sigma
\\
&\underset{\T{\eqref{increment} (ii)}}\simleq\int_0^1\Big|\di_\sigma F^*_\alpha(\sigma e^{i\theta})\Big|d\sigma
\\
&\underset{\T{\eqref{increment} (i)}}=\Big|\int_0^1\di_\sigma F^*_\alpha(\sigma(1-e^{i\theta}))d\sigma\Big|
\\
&=\Big|F^*_\alpha(1-e^{i\theta})\Big|.
\end{split}
\end{equation}
It follows
\begin{equation*}
\begin{split}
\Big|h\Big(\Big(F^*_\alpha(1-e^{i\theta}+\frac1\nu)-F^*_\alpha(\frac1\nu)\Big),0,r \Big)\Big|&\underset{\T{\eqref{1.1}}}=O(F(| F^*_\alpha(1-e^{i\theta})|))
\\
&\simleq |\theta|^\alpha.
\end{split}
\end{equation*}
Thus $\frac{|h(F^*_{\alpha\,\nu},0,r )|}{1-\cos\theta}\simleq |\theta|^{-2+\alpha}=|\theta|^{-1+\beta},\,\,\beta>0$; we can therefore apply the dominated convergence theorem to the sequence $\{\di_tv_\nu\}_\nu$ and conclude that $\di_tv_\nu|_{t=1}\to\di_tv|_{t=1}$. It follows that for $\nu$ large,  the disc over $F^*_{\alpha\,\nu}(1-\Delta)$ is transversal to $b\Om$, that is
$$
\di_tv_\nu>0.
$$
Thus the disc $\mathcal S_\nu$ attached to $b\Om$ over the sector $F^*_{\alpha\,\nu}(1-\Delta)$ ``points down" at $\tau=1$ for $\nu$ large. 
By the aid of the discs attached to $b\Om$ over a family of translations of the sector $F^*_{\alpha\,\nu}(1-\Delta)$, 
we sweep out a full neighborhood of $0$ in the complement of $\Om$ and thus get the extension at $0$ of a holomorphic function (if this extends at $F^*_\alpha(2)$).

\epf
\br
\Label{r1.1}
If $h|_{S_\alpha\times\{0\}\times(-\epsilon,\epsilon)}\le0$, that is, $S_\alpha\times\{0\}\times(-\epsilon.\epsilon)\subset\bar\Om$, then the component $v$ of the disc $\mathcal S_\alpha$ satisfies $\di_tv\ge0$; if, moreover $h<0$ at some point of $b\mathcal S_\alpha$, then $\di_tv>0$.
Hence we have in hands from the beginning a disc satisfying \eqref{1.5}. This simplifies the proof of Theorem~\ref{t1.1} and also dispense from taking the bumping $h-\eta\chi$. In other words, we have forced extension at 0 of a holomorphic function $f$ on $\Om$: in order that a holomorphic function $f$ extends at $0$, it needs not to extend at some other point of $b\mathcal S$.
\er

\section{Pseudoconvex bumps}
\Label{s2}
Let $F^*$ be a holomorphic, injective, function of the $\epsilon$-neighborhood of $0$ in a sector of $\C$ with axis $\R^+$ and angle $>\pi$ whose range contains $\R^+$ at $0$,   denote by $F$ the inverse to $F^*$, and let $\Om$ be a domain in $\C^n$ defined by $s>h(z,z',r )$ with $h=O(F)$. We also make the assumptions \eqref{1.1,5} and \eqref{1.1,6} and \eqref{1.1,7}  on $h$ and $F^*$ which yield $h_r(F^*_\alpha)\in C^{1,\beta}$ for $h_r=h(\cdot,0,r)$ according to Lemma~\ref{l1.1}.
 We set $F^*_\alpha(z):=F^*(\epsilon z^\alpha)$ and consider  the sector $S_\alpha=\{z:\,z=F^*_\alpha(1-\tau),\,\tau\in\Delta\}$ and the disc $\mathcal S_\alpha $ attached to $b\Om$ over $S_\alpha\times\{0\}$, that is, $\mathcal S_\alpha=\{(F^*_\alpha(1-\tau),0,u(\tau)+iv(\tau)):\,\tau\in\Delta\}$ where $u$ and $v$ are the harmonic extension of the functions satisfying $v-h(F^*_\alpha(1-\tau),0,u)|_{b\Delta}=0$ and $v|_{b\Delta}=T_1u|_{b\Delta}$. We will suppose, all through this section, that $\alpha<1$. We take $\alpha_1$ with $\alpha<\alpha_1<1$; the crucial and elementary remark which underlies this part of the discussion is that, since $F\circ F^*$ is  the identy of $\C^+$, then $\Re F>0$ on $\C^+$ and moreover, since $\Re z\sim |z| $ for $z=\epsilon(1-\tau)^{\alpha_1}$, then 
\begin{equation}
\Label{2.1bis}
\Re F\sim|F|\quad\T{on $S_{\alpha_{1}}$}.
\end{equation}
\bp
\Label{p2.1}
Let $h=h(z)$ with $h(0)=0$ be plurisubharmonic, let $F$ be holomorphic and assume, uniformly on $z'$ and $r$
\begin{equation}
\Label{2.2}
\di\dib h\simgeq \frac{|\di_zF|^2}{|F|}\qquad\T{on $S_{\alpha_{1}}\setminus S_{\alpha}$ for $\alpha<\alpha_1<1$}.
\end{equation}
Then, there is $\tilde h$ with $\tilde h(0)=0$, which is subharmonic and satisfies
\begin{equation}
\Label{2.3}
\begin{cases}
\tilde h\leq h,
\\
\tilde h<h\quad\T{on $S_{\alpha}$.}
\end{cases}
\end{equation}
\ep
\bpf
We define 
$$
\tilde h=h-\eta \chi_{S_{\alpha_{1}}\setminus S_{\alpha}}\Re F,
$$
where the cut-off $\chi_{S_{\alpha_{1}}\setminus S_{\alpha}}$ is the pull-back under $F$ of the conical cut-off  which is $0$ for $|\arg z|\ge \frac{\alpha_1\pi}2$ and $1$ for $|\arg z|\le \frac{\alpha\pi}2$,
that is, $\chi_{S_{\alpha_{1}}\setminus S_{\alpha}}=\chi(\frac{F-\bar F}{F+\bar F})$. Thus
\begin{equation}
\begin{split}
\Label{2.4}
\di_z \chi\Big(\frac{F-\bar F}{F+\bar F}\Big)&=\dot \chi\frac{\bar F\di_zF}{(F+\bar F)^2}
\\
&\simleq \frac{|\di_zF|}{|F|},
\end{split}
\end{equation}
and
\begin{equation}
\Label{2.5}
\begin{split}
\Big|\di_{\bar z}\di_{ z}\chi\Big(\frac{F-\bar F}{F+\bar F}\Big)\Big|&\simleq \Big|\frac{\bar F\di_zF}{F+\bar F}\Big|^2+\Big|\frac{|\di_zF|^2(F+\bar F)^2-|\di_zF|^2\bar F(F+\bar F)^2}{(F+\bar F)^4}\Big|
\\
&\simleq \frac{|\di_zF|^2}{|F|^2}.
\end{split}
\end{equation}
It then follows
\begin{equation}
\Label{2.6}
\begin{split}
\di_z\di_{\bar z}(\chi_{S_{\alpha_{1}}\setminus S_{\alpha}}\Re F)&=\di_z\di_{\bar z}\chi_{S_{\alpha_{1}}\setminus S_{\alpha}}\Re F-2\Re( \di_z\chi_{S_{\alpha_{1}}\setminus S_{\alpha}})\di_z\Re F
\\
&\underset{\T{\eqref{2.4} and \eqref{2.5}}}\simleq \frac{|\di_zF|^2}{|F|}
\\
&\underset{\T{\eqref{2.2}}}\simleq \di_z\di_{\bar z} h\quad\T{on supp$(\chi_{S_{\alpha_{1}}\setminus S_{\alpha}})$.}
\end{split}
\end{equation}
From \eqref{2.6} we readily conclude that $\di_z\di_{\bar z} \tilde h\geq0$; thus $\tilde h$ is subharmonic. It also satisfies the other requirements of the statement.

\epf
As an immediate consequence of Proposition~\ref{p2.1} we have the proof of
\bt
\Label{t2.1}
Let $\Om$ be a pseudoconvex domain defined by $s=h(z,z',r )$ with $h$ satisfying \eqref{2.2}, let $\tilde h$ be obtained by the technique of Proposition~\ref{p2.1}, and let 
$\tilde \Om$ be defined by $s>\tilde h(z,z',r )$. Then $\tilde \Om$ has the following properties
\begin{equation*}
\begin{cases}
\tilde\Om\T{ is pseudoconvex},
\\
\tilde \Om\supset \Om\cup(b\Om\cap\pi^{-1}(S_{\alpha}\times\{0\}\times\{0\})\quad\T{where $\pi$ is the projection to $w=0$},
\\
0\notin\tilde \Om.
\end{cases}
\end{equation*}
\et
In particular, there is a function which is holomorphic in $\Om$, extends holomorphically across $b\Om$ at any point of the boundary of the attached disc $b\mathcal S_{\alpha}\setminus\{0\}$, but is singular at $0$; thus the sector $\mathcal S_{\alpha}$ attached to $b\Om$ over $S_\alpha\times\{0\}\times\{0\}$, even in case  is contained in $\Om$ or is tangent to $b\Om$  at $0$, nonetheless  is not a propagator of holomorphic extendibility at the vertex $0$.

\end{document}